\input amstex
\magnification=\magstep1 
\baselineskip=13pt
\documentstyle{amsppt}
\vsize=8.7truein \CenteredTagsOnSplits \NoRunningHeads
\def\zz{\bold{z}}
\def\EE{\bold {E\thinspace}}
\def\dist{\operatorname{dist}}
\def\supp{\operatorname{supp}}

\def\ind{\operatorname{ind}}

\def\LL{\Cal{L}}
\def\PP{\Cal{P}}
\def\TT{\Bbb{T}}
\topmatter
 
 \title Weighted counting of solutions to sparse systems of equations \endtitle 
\author Alexander Barvinok and Guus Regts \endauthor
\address Department of Mathematics, University of Michigan, Ann Arbor,
MI 48109-1043, USA \endaddress
\email barvinok$\@$umich.edu \endemail
\address Korteweg de Vries Institute for Mathematics, University of Amsterdam \endaddress
\email guusregts$\@$gmail.com \endemail
\date February  2019 \enddate
\thanks  The research of the first author was partially supported by NSF Grant DMS 1361541. The research of the second author was supported by 
a personal NWO Veni grant.
\endthanks 
\abstract  Given complex numbers $w_1, \ldots, w_n$, we define the weight $w(X)$ of a set $X$ of 0-1 vectors as the sum of 
$w_1^{x_1} \cdots w_n^{x_n}$ over all vectors $(x_1, \ldots, x_n)$ in $X$. We present an algorithm, which for a set $X$ defined by a system of homogeneous linear equations with at most $r$ variables per equation and at most $c$ equations per variable, computes $w(X)$ within relative error $\epsilon >0$ in $(rc)^{O(\ln n-\ln \epsilon)}$ time provided $|w_j| \leq \beta (r \sqrt{c})^{-1}$ for an absolute constant $\beta >0$ and all $j=1, \ldots, n$. A similar algorithm is constructed for computing the weight of a linear code over ${\Bbb F}_p$.  Applications include counting weighted perfect matchings in hypergraphs, counting weighted graph homomorphisms, computing weight enumerators of linear codes with sparse code generating matrices, and computing the partition functions of the ferromagnetic Potts model at low temperatures and of the hard-core model at high fugacity on biregular bipartite graphs.
\endabstract
\endtopmatter
\document

\head 1. Weighted counting of 0-1 vectors\endhead

\subhead (1.1) Weight of a set of 0-1 vectors \endsubhead
Let us  fix complex numbers $w_1, \ldots, w_n$, referred to as {\it weights} in what follows.
We define the {\it weight} $w(x)$ of a 0-1 vector $x \in \{0, 1\}^n$ by 
$$w(x)=w_1^{\xi_1} \cdots w_n^{\xi_n} =\prod_{j:\ \xi_j=1} w_j \quad \text{where} \quad x=\left(\xi_1, \ldots, \xi_n \right).$$
Here we agree that $0^0=1$, so that $w(x)$ is a continuous function of $w_1, \ldots, w_n$ for a fixed $x$.

We define the {\it weight} of a finite set $X \subset \{0, 1\}^n$ by 
$$w(X) = \sum_{x \in X} w(x) = \sum\Sb x \in X, \\ x=\left(\xi_1, \ldots, \xi_n\right) \endSb w_1^{\xi_1} \cdots w_n^{\xi_n}. \tag1.1.1$$

Given $X \subset \{0, 1\}^n$, the value of $w(X)$ as a function of $w_1, \ldots, w_n$ is also known as the {\it partition function}  or {\it generating function} of $X$.

Our first main result is as follows.

\proclaim{(1.2) Theorem} Let $A=\left(a_{ij}\right)$ be an $m \times n$ integer matrix and let us define
$X \subset \{0, 1\}^n$ by 
$$ X=\Bigl\{x \in \{0, 1\}^n,\ x=(\xi_1, \ldots, \xi_n):\quad \sum_{j=1}^n a_{ij} \xi_j = 0 \quad \text{for} \quad i=1, \ldots, m \Bigr\}.$$
Suppose that the number of non-zero entries in every row of $A$ does not exceed $r$ for some $r \geq 2$ and that the number of non-zero entries in every column of $A$ does not exceed $c$ for some $c \geq 1$. There is an absolute constant $\alpha >0$ such that if $w_1, \ldots, w_n \in {\Bbb C}$ are weights satisfying
$$\left| w_j \right| \ \leq \ {\alpha \over r \sqrt{c}} \quad \text{for} \quad j=1, \ldots, n,$$
then
$$w(X) \ne 0.$$ One can choose $\alpha=0.46$.
\endproclaim

Geometrically, the set $X$ in Theorem 1.2 is the set of 0-1 vectors in a subspace. We are interested in efficient algorithms to compute $w(X)$ approximately. Theorem 1.2 implies that such an efficient algorithm exists for a non-trivial range of weights $w_1, \ldots, w_n$ provided the matrix $A$ is sufficiently sparse (that is, $r$ and $c$ are sufficiently small), even when the dimension $n$ of the ambient space is allowed to be large. This connection between the sparsity condition for $A$ (frequent in applications and easily verified) and the computational complexity of $w(X)$ appears to be new.

\subhead (1.3) Computing $w(X)$  \endsubhead Theorem 1.2 implies that $w(X)$ can be efficiently approximated as long as the weights $w_j$ satisfy a slightly stronger inequality, 
$$\left| w_j \right| \ \leq \ {\beta \over r \sqrt{c}} \quad \text{for} \quad j=1, \ldots, n\tag1.3.1$$
for any $\beta < \alpha$, fixed in advance, so one can choose $\beta=0.45$.
We describe the connection below, see also Section 1.2 of \cite{Ba16}.

Without loss of generality we assume that the matrix $A$ has no zero rows and no zero columns (although this assumption is not needed in this section, it will be relevant later in Section 5). Indeed, zero rows of $A$ can be ignored and if, say, the $n$-th column of $A$ is zero, we have 
$$w(X)=\left(1+w_n \right) w\left(\widehat{X}\right),$$
where $\widehat{X} \subset \{0, 1\}^{n-1}_+$ is the set defined by the system $\widehat{A} x=0$, where $\widehat{A}$ is the $m \times (n-1)$ matrix obtained from $A$ by deleting the $n$-th column. 

 For a $\zeta \in {\Bbb C}$, let $\zeta w_1, \ldots, \zeta w_n$ be the scaling of the weights and let 
$w(X; \zeta)$ be the corresponding weight of $X$ so that $w(X; 1)=w(X)$ while $w(X; 0)=1$ (note that $0 \in X$). Theorem 1.2 implies that as long as the weights $w_j$ satisfy (1.3.1), we have
$$w(X; \zeta) \ne 0 \quad \text{provided} \quad |\zeta| \ \leq \ {\alpha \over \beta} =:\gamma. \tag1.3.2$$ 
Note that $\gamma >1$.

Let us choose a continuous branch of $f(\zeta)=\ln w(X; \zeta)$ for $|\zeta| \leq \gamma$ and let 
$$T_s(\zeta) = f(0) + \sum_{k=1}^s {f^{(k)}(0) \over k!} \zeta^k \tag1.3.3$$
be the Taylor polynomial of $f$ of some degree $s$ computed at $\zeta=0$. Since (1.3.2) holds and $w(X; \zeta)$ is a polynomial of degree at most $n$ in 
$\zeta$, we have 
$$\left| f(1) - T_s(1) \right| \ \leq \ {n \over (s+1) \gamma^s (\gamma-1)}, $$
see Lemma 2.2.1 of \cite{Ba16}. Using that $\gamma>1$, we conclude that to approximate $f(1) =\ln w(X)$ within an additive error $\epsilon >0$ by $T_s(1)$, it suffices to choose $s=O\left(\ln n - \ln \epsilon\right)$, where the implied constant in the ``$O$" notation depends only on $\gamma$. We say then that 
$e^{T_s(1)}$ approximates $w(X)$ within relative error $\epsilon$.

We have $f(0)=0$ and computing $f^{(k)}(0)$ for $k=1, \ldots, s$ reduces to computing 
$${d^k \over d \zeta^k} w(X; \zeta) \Big|_{\zeta=0} \quad \text{for} \quad k=1, \ldots, s \tag1.3.4$$
in $O(s^2)$ time. Indeed, it is not hard to see 
that the values $f^{(k)}(0)$ are the solutions of a non-degenerate triangular system of linear equations with right hand side given by (1.3.4),  
see Section 2.2.2 of \cite{Ba16}.
Furthermore,
$${d^k \over d \zeta^k} w(X; \zeta) \Big|_{\zeta=0} = k! \sum\Sb x \in X, \\ x=\left(\xi_1, \ldots, \xi_n\right): \\ \xi_1 + \ldots + \xi_n =k \endSb 
w_1^{\xi_1} \cdots w_n^{\xi_n},$$
so computing (1.3.4) reduces to the inspection of all points $x \in X$, $x=\left(\xi_1, \ldots, \xi_n\right)$, satisfying $\xi_1 + \ldots + \xi_n \leq s$, which can be 
done through the exhaustive search in $mn^{O(s)}$ time. Given that $s=O(\ln n - \ln \epsilon)$, this produces an algorithm approximating $w(X)$ within a 
relative error $\epsilon >0$ in quasi-polynomial $n^{O(\ln n - \ln \epsilon)}$ time, where the implied constant in the ``$O$" notation depends only on $\gamma$ in (1.3.2). In Section 5 we show that we can 
compute $f^{(k)}(0)$ in (1.3.3) faster, in $(rc)^{O(\ln n-\ln \epsilon)}$ time. In particular, if $r$ and $c$ are fixed in advance, we obtain a polynomial time approximation algorithm.

Next, we consider enumerating 0-1 vectors in affine subspaces, not necessarily containing the origin.

\subhead (1.4) Non-homogeneous linear equations in 0-1 vectors \endsubhead  
We interpret a vector $x=\left(\xi_1, \ldots, \xi_n\right)$ as a column $n$-vector. Let $A$ be an $m \times n$ integer matrix as above, let $b$ be an integer $m$-vector and let 
$$X=\Bigl\{x \in \{0, 1\}^n:\quad Ax =b \Bigr\}$$
be the set of 0-1 vectors satisfying a system of linear equations with matrix $A$.
In general, it is an NP-hard problem to decide whether $X$ is empty, so there is no hope to compute $w(X)$ efficiently. 

Suppose, however, that we are presented with a point $y \in X$, $y=\left(\eta_1, \ldots, \eta_n \right)$. 
Every point $x \in X$ can be uniquely written as $x=y +z$, $z=\left(\zeta_1, \ldots, \zeta_n\right)$,
where $Az=0$ and $\zeta_j \in \{-1, 0\}$ if $\eta_j=1$ and $\zeta_j \in \{0, 1\}$ if $\eta_j =0$. Let $a_1, \ldots, a_n$ be the columns of $A$ and let 
$\widehat{A}$ be the matrix obtained from $A$ by replacing $a_j$ with $-a_j$ whenever $\eta_j=1$. Let 
$$Z =\Bigl\{ z \in \{0, 1\}^n: \quad \widehat{A}z =0\Bigr\}.$$
Hence every point $x \in X$, $x=\left(\xi_1, \ldots, \xi_n\right)$,
can be uniquely written as $\xi_j = \eta_j + \sigma_j \zeta_j$, where for $z=\left(\zeta_1, \ldots, \zeta_n\right)$ we have $z \in Z$ and 
$$\sigma_j=\cases 1 &\text{if\ } \eta_j=0, \\ -1 &\text{if\ } \eta_j =1. \endcases$$
Then, for the weight of $Z$, we have 
$$
\aligned w(Z) = &\sum\Sb z \in Z, \\ z=\left(\zeta_1, \ldots, \zeta_n\right) \endSb \prod_{j=1}^n w_j^{\zeta_j} = \sum\Sb x \in X, \\ x=\left(\xi_1, \ldots, \xi_n\right) \endSb \prod_{j=1}^n w_j^{\sigma_j(\xi_j -\eta_j)} 
\\
=&\sum\Sb x \in X, \\ x=\left(\xi_1, \ldots, \xi_n \right) \endSb \prod_{j: \xi_j \ne \eta_j} w_j.\endaligned \tag1.4.1$$
For $x \in \{0, 1\}^n$, 
$x=\left(\xi_1, \ldots, \xi_n\right)$, let 
$$\dist(x, y) = \left| \left\{j:\ \xi_j \ne \eta_j \right\} \right|  $$
be the Hamming distance between $x$ and $y$.

In particular, if we choose 
$$w_1 = \ldots = w_n =\omega$$
for some $\omega$, we get 
$$w(Z)= \sum_{x \in X} \omega^{\dist(x, y)}. \tag1.4.2$$
Assuming that every row of $A$ contains not more than $r \geq 2$ non-zero entries and every column of $A$ contains not more than $c \geq 1$ non-zero entries, we conclude that the sum (1.4.2) can be computed within relative error $\epsilon >0$ in $(rc)^{O(\ln n-\ln \epsilon)}$ time provided
$$|\omega| \ \leq \ {\beta \over r \sqrt{c}},$$
where $\beta >0$ is an absolute constant (one can choose $\beta=0.45$). If $r$ and $c$ are fixed in advance, we have a polynomial time 
approximation algorithm of $(n/\epsilon)^{O(1)}$ complexity.

In the next section we consider combinatorial applications of our result. 
We first consider a variation of Theorem 1.2 that applies to codes.

\subhead (1.5) Weight of a code \endsubhead
Let $\kappa > 1$ be an integer. We consider $n$-vectors $x=\left(\xi_1, \ldots, \xi_n\right)$ with coordinates 
$\xi_j$ taking values in the set $\{0, \ldots, \kappa -1\}$, which we interpret as the set ${\Bbb Z}/\kappa {\Bbb Z}$ of remainders modulo $\kappa$. Given $n$ complex numbers $w_1, \ldots, w_n$, we 
define the {\it weight} $w(x)$ of a vector $x \in \left({\Bbb Z}/\kappa {\Bbb Z}\right)^n$ by 
$$w(x)=\prod_{j:\ \xi_j \ne 0} w_j \quad \text{for} \quad x=\left(\xi_1, \ldots, \xi_n\right)$$ 
and the {\it weight} $w(X)$ of a set $X \subset \left({\Bbb Z}/\kappa {\Bbb Z}\right)^n$ by 
$$w(X)=\sum_{x \in X} w(x)$$
(we agree that the weight of the zero vector is 1).

We obtain the following result.
\proclaim{(1.6) Theorem} Let $A=\left(a_{ij}\right)$ be an $m \times n$ integer matrix and let us define a set $X \subset \left({\Bbb Z}/\kappa {\Bbb Z}\right)^n$ by 
$$\split X=\Bigl\{x \in \left({\Bbb Z}/\kappa {\Bbb Z}\right)^n,\ x=\left(\xi_1, \ldots, \xi_n\right): \quad &\sum_{j=1}^n a_{ij} \xi_j \equiv 0 \mod \kappa \\
&\quad \text{for} \quad i=1, \ldots, m\Bigr\}. \endsplit$$
Suppose that the number of non-zero entries in every row of $A$ does not exceed $r$ for some $r \geq 2$ and that the number of non-zero entries in every column of $A$ does not exceed $c$ for some $c \geq 1$. There is an absolute constant $\alpha >0$ such that if $w_1, \ldots, w_n \in {\Bbb C}$ are weights satisfying 
$$\left| w_j \right| \ \leq \  {\alpha \over (\kappa-1) r \sqrt{c}} \quad \text{for} \quad j=1, \ldots, n,$$
then 
$$w(X) \ne 0.$$
One can choose $\alpha =0.46$.
\endproclaim
 As in Section 1.3, we obtain an algorithm of $(rc)^{O(\ln \kappa n -\ln \epsilon)}$ complexity to approximate $w(X)$ within relative error $\epsilon >0$ provided 
$$|w_j| \ \leq \ {\beta \over (\kappa-1) r \sqrt{c}} \quad \text{for} \quad j=1, \ldots, n,$$
where $\beta < \alpha$ is fixed in advance (we can choose $\beta=0.45$). For $r$ and $c$ fixed in advance, the algorithm has polynomial $(\kappa n/\epsilon)^{O(1)}$ complexity. 

\subsubhead Organization \endsubsubhead
We deduce Theorem 1.2 and Theorem 1.6 from a general result asserting that 
$$\int_{\TT^m} e^{p(z)} \ d \mu \ne 0,$$
for some Laurent polynomials $p: \TT^m \longrightarrow {\Bbb C}$ on the torus $\TT^m$ endowed with a product probability measure $\mu$
(see Theorem 3.2 and Corollary 3.3 below).
 After that, the proofs of Theorems 1.2 and 1.6 are completed in a more or less straightforward way in Section 4.
 
 In Section 5, we provide details of an approximation algorithm for $w(X)$.
 We do not discuss an analogous algorithm for codes in Theorem 1.6 as it is very similar.
We first consider some concrete combinatorial applications of these results in Section 2 below.

\head 2. Combinatorial applications \endhead

We apply Theorem 1.2 to weighted counting of perfect matchings in hypergraphs, computing the partition function of the hard-core model at high fugacity for biregular bipartite graphs and  to weighted counting of graph homomorphisms. We apply Theorem 1.6 to computing weight enumerators of linear codes with sparse code generating matrices and to computing the partition function of the ferromagnetic Potts model at low temperatures.

\subhead (2.1) Perfect matchings in hypergraphs \endsubhead A {\it hypergraph} $H=(V, E)$ is a finite set $V$ of {\it vertices} together with a collection $E$ of non-empty subsets $V$, called {\it edges} of the hypergraph. The {\it degree} of a vertex $v$ is the number of edges $e \in E$ that contain $v$.
A {\it perfect matching} in $H$ is a set of pairwise disjoint edges $e_1, \ldots, e_n$, such that $e_1 \cup \ldots \cup e_n = V$. Let us introduce a 0-1 variable 
$x_e$ for each $e \in H$. We encode a collection of edges of $H$ by a 0-1 vector, where
$$x_e=\cases 1 &\text{if $e$ is in the collection} \\ 0 &\text{otherwise.} \endcases$$
Then $e_1, \ldots, e_n$ is a perfect matching if and only if
$$\sum_{e:\ v \in e} x_e =1 \quad \text{for all} \quad v \in V. \tag2.1.1$$
In the system (2.1.1) the number of variables per equation is the maximum degree $d$ of a vertex of $H$ and the number of equations per variable is the maximum cardinality $k$ of an edge. It is an NP-complete problem to find if a given hypergraph contains a perfect matching provided $k \geq 3$, see, for example, Problem SP1 in \cite{A+99}. However, as follows from 
Section 1.4, given {\it one} perfect matching $M_0$, we can efficiently approximate a certain statistic over {\it all} perfect matchings $M$ of $H$, namely 
the sum 
$$\sum_{M\in {\Cal M}(H)} \omega^{\dist(M_0, M)}, \tag2.1.2$$
where ${\Cal M}(H)$ is the set of all perfect matchings, $\dist(M_0, M)$ is the Hamming distance 
between matchings, that is, the number of edges where the matchings differ and
$$|\omega| \ \leq \ {\beta \over d\sqrt{k}}.$$
The complexity of the algorithm approximating 
(2.1.2) within relative error $\epsilon >0$ is $(d k)^{O(\ln |E|-\ln \epsilon)}$. If $d$ and $k$ are fixed in advance, the algorithm 
achieves polynomial $(|E|/\epsilon)^{O(1)}$ complexity.
This can be contrasted with the fact that knowing one solution of a problem generally does not help to find another or to count all solutions, cf. \cite{Va79} and \cite{VV86}.

Is is shown in \cite{Ba18} that if the hypergraph is {\it uniform} and $k$-{\it partite}, that is, we have $V=V_1 \cup \ldots \cup V_k$ with pairwise disjoint $V_1, \ldots, V_k$ 
such that $|V_1| = \ldots = |V_k|=n$ and every edge $e \in E$ contains exactly one vertex from each $V_i$, then one can efficiently approximate (2.1.2) under the weaker condition 
$$|\omega| \ \leq \ {\beta \over \sqrt{d-1}}$$
for any $\beta < 1$, fixed in advance.

\subhead (2.2) The hard-core model at high fugacity \endsubhead
Given an undirected graph $G=(V, E)$, a set $S \subset V$ of vertices is called {\it independent} if no two vertices of $S$ span an edge of $G$ (we agree that $S=\emptyset$ is 
always independent). The {\it independence polynomial} of $G$ is a univariate polynomial defined by
$$p_G(\lambda) =\sum\Sb S \subset V \\ S \text{\ is independent} \endSb \lambda^{|S|}, \tag2.2.1$$
see for example, Chapter 6 of \cite{Ba16}. It is also known the partition function of the hard-core model. The parameter $\lambda$ is known as the {\it fugacity}.

The problem of (approximately) computing the number of independent sets in a bipartite graph is considered to be computationally hard. It is the basis of the class of $\#$BIS hard problems, and it is known that to approximate the $p_G(\lambda)$ on bipartite graphs of maximum degree $d$ is a $\#$BIS hard problem, provided  $\lambda>\frac{(d-1)^{d-1}}{(d-2)^d}$~\cite{C+16}. Moreover, as the authors of \cite{C+16} informed us, it follows form their construction that computing $p_G(\lambda)$ for sufficiently large $\lambda$ remains a $\#$BIS-hard problem when restricted to bipartite $d$-regular graphs $G$. 

In \cite{J+19} it was however shown that for $d\geq 3$ there exists $\lambda^{\ast}=\lambda^{\ast}(d)>0$ such that for all $\lambda>\lambda^{\ast}$ and all $d$-regular, bipartite, expander graphs $G$, the value of $p_G(\lambda)$ can be approximated in polynomial time.
Here we will use Theorem 1.2 to show that for each fixed $d_1,d_2\in \Bbb{N}$ such that $d_2-d_1\geq 1$, there exists $\lambda_0=\lambda_0(d_1, d_2)>0$ such that for all $\lambda>\lambda_0$ and any biregular, bipartite graph with degrees $d_1,d_2$ we can approximate $p_G(\lambda)$ in polynomial time.

To this end, let us fix a biregular bipartite graph $G=(V,E)$ with degrees $d_1$ and $d_2\geq d_1+1$. We write $V=L\cup R$ for the bipartition and we assume that each vertex in $L$ has degree $d_1$ and each vertex in $R$ has degree $d_2$.
For an independent set $I$ we write $I_L:=I\cap L$ and $I_R:=I\cap R$.

We wish to encode $p_G(\lambda)$ as the weight $w(X)$ of a suitably defined set $X$. We direct all edges from $L$ to $R$, thus making $G$ a directed graph.
We associate to each vertex $v\in V$ a 0-1 variable $x_v$ and to each edge $(u, v) \in E$ a 0-1 variable $x_{uv}$.
Let $X$ be the solution set to the following system of equations:
$$ -x_u+x_v+x_{uv}=0 \quad \text{ for each directed edge} \quad  (u,v)\in E. \tag2.2.2$$
Any $x\in X$ uniquely corresponds to an independent set $I$ of $G$. 
Indeed, let $I$ be the the sets of vertices $u\in L$ for which $x_u=0$ and vertices $v\in R$ for which $x_v=1$. 
Then for $u\in I_L$ none of its neighbors will be contained in $I$ since for each edge $(u,v)$, the value of $x_v$ is forced to be zero.
Similarly, for any $v\in I_R$, none of its neighbors will be contained in $I$ since for each edge $(u,v)$, the value of $x_u$ is forced to be $1$.
Hence the set $I$ is independent. 
Conversely, if $I$ is an independent set, setting  
$$\split &x_u=\cases 0 &\text{if\ } u \in I_L, \\ 1 &\text{if\ } u \in L \setminus I_L, \endcases  \quad x_v=\cases 1 &\text{if\ } v \in I_R, \\ 0 &\text{if\ } v \in R \setminus I_R \endcases 
\quad \text{and} \\
 &x_{uv}=\cases 0 &\text{if\ } u \in I_L  \quad \text{or} \quad v \in I_R, \\ 1 &\text{if\ } u \in L \setminus I_L \quad \text{and} \quad v \in R \setminus I_R, \endcases \endsplit$$
 gives a solution to (2.2.2).

Next, we introduce weights $w_u$ for the coordinates $x_u$ with $u\in L$, weights $w_v$ for the coordinates $x_v$ with $v \in R$ and weights $w_{uv}$ for the coordinates $x_{uv}$ with $(u, v) \in E$ as follows:
$$\split &w_u = \omega^{(d_2-d_1)/2} \quad \text{for} \quad u \in L, \\ &w_v=\omega^{(d_2-d_1)/2} \quad \text{for} \quad v \in R  \quad \text{and} \\ &w_{uv}=\omega \quad \text{for} \quad 
(u, v) \in E. \endsplit$$
For a solution $x\in X$ corresponding to an independent set $I$, we then have
$$
\split
w(x)&=\left(\prod_{v\in L\setminus I_L}\omega^{(d_2-d_1)/2}\right) \left( \prod_{\Sb\{u,v\}\in E\\ u,v\notin I\endSb} \omega\right) \left( \prod_{u\in I_R}\omega^{(d_2-d_1)/2} \right)
\\
&= \omega^{(d_2-d_1)(|L|-|I_L|)/2} \cdot \omega^{|E|-d_1|I_L|-d_2|I_R|} \cdot \omega^{(d_2-d_1)|I_R|/2} 
\\
&=\omega^{(d_2-d_1)|L|/2 +|E|}\cdot \omega^{-(d_1+d_2)|I_L|/2}  \cdot \omega^{-(d_1+d_2)|I_R|/2} \\
&=\omega^{(d_1+d_2)|L|/2} \omega^{-(d_1+d_2) |I|/2}.
\endsplit
$$
In other words, for the weight of $X$, we have 
$$w(X)=\omega^{(d_1+d_2)|L|/2} p_G\left(\frac{1}{\omega^{(d_1+d_2)/2}}\right)$$
for the independence polynomial $p_G$ defined by (2.2.1).

Now, since in (2.2.2) the number of variables per equation is $3$ and the number of equations per variable is at most $d_2$, it follows from Theorem 1.2 that if 
$$|\lambda| \ \geq \  \left(6.7 \sqrt{d_2}\right)^{d_1+d_2} \ > \ \left(\frac{3\sqrt{d_2}}{0.45}\right)^{d_1+d_2},$$ 
then $p_G(\lambda)\neq 0$ and moreover that we can efficiently approximate $p_G$ (in polynomial time if $d_2$ is fixed in advance).
We moreover note that with a similar argument, for a $d$-regular bipartite graph $G=(L\cup R,E)$, we can efficiently approximate the sum $$\sum_{\Sb I \subset L\cup R \\ I \text{ is independent} \endSb}\lambda^{| I\cap L |}$$ for large $\lambda.$
This is somewhat similar in spirit to a result of van den Berg and Steiff~\cite{BS94}, who showed that for the integer lattice $\Bbb{Z}^d$, assigning $\lambda_1>0$ to vertices with even coordinate sum and $\lambda_2>0$ to vertices with odd coordinate sum, for all but a countable set of pairs $(\lambda_1,\lambda_2)$ the associated Gibbs measure is unique.

\subhead (2.3) Weighted counting of graph homomorphisms \endsubhead
Let $G_1=(V_1, E_1)$ be an undirected graph without loops or multiple edges and let $G_2=(V_2, E_2)$ be an undirected graph without multiple edges, but possibly with loops. We assume that $V_2=\{1, \ldots, n\}$ and assume that $G_1$ and $G_2$ are both connected. 
 A map $\phi: V_1 \longrightarrow V_2$ is called a {\it homomorphism} if $\phi(u)$ and $\phi(v)$ span an edge of $G_2$ whenever $u$ and $v$ span an edge of $V_1$. If $V_2$ is the complete graph without loops then every homomorphism $\phi: G_1 \longrightarrow G_2$ is naturally interpreted as a coloring of 
 the vertices of $G_1$ with a set of $n$ colors such that no two vertices spanning an edge of $G_1$ are colored with the same color (such colorings are called 
 {\it proper}). For any fixed $n \geq 3$, it is an NP-complete problem to decide wether a given graph admits a proper $n$-coloring, see for example, 
 Problem GT5 in \cite{A+99}.
Our goal is to encode all homomorphisms $\phi: G_1 \longrightarrow G_2$ that map a fixed vertex $a \in V_1$ to a fixed vertex, say $n$, of $G_2$ as the set of 0-1 solutions to a system of linear equations.

We say that vertices $u, v \in V_1$ are {\it neighbors} if $\{u, v\} \in E_1$.
We orient the edges of $G_1$ arbitrarily, so that an edge of $G_1$ is an ordered pair of neighbors $(u, v)$.  Let us introduce 0-1 variables $x^{uv}_{ij}$ indexed by (now directed) edges $(u, v) \in E_1$ and ordered pairs $1 \leq i, j \leq n$ such that $\{i, j\} \in E_2$ (we may have $i=j$). 
The idea is to use the variables $x^{uv}_{ij}$ to encode a map $\phi: V_1 \longrightarrow V_2$, so that 
$$x^{uv}_{ij}=\cases 1 &\text{if\ } \phi(u)=i \quad \text{and} \quad \phi(v)=j, \\ 0 &\text{otherwise.} \endcases
\tag2.3.1$$
For every ordered pair of neighbors $(u, v)$ and every vertex $i \in V_2$ we define the sum
$$\split S^{u,v}_i =&\sum_{j:\ \{i, j\} \in E_2}  x^{uv}_{ij} \quad \text{if} \quad (u, v ) \in E_1 \quad \text{and} \\
S^{u,v}_i= &\sum_{j:\ \{i, j\} \in E_2} x^{vu}_{ji} \quad \text{if} \quad (v, u) \in E_1 \endsplit \tag2.3.2$$
and for every $u \in V_1$ and every $i \in V_2$, we introduce the following equations:
\bigskip
\noindent (2.3.3) Fix $u \in V_1 \setminus \{a\}$ and $i \in V_2$. The sums $S^{u,v}_i$, where $v$ is a neighbor of $u$, are all equal.
\bigskip
\noindent The idea, of course, is that the sums (2.3.2) are all equal to $1$ if $\phi(u)=i$ and equal to $0$ if $\phi(u) \ne i$.
Next, we encode the condition $\phi(a)=n$ by the following system of equations:
\bigskip
\noindent (2.3.4) For all neighbors $v$ of $a$,  
$$S^{a,v}_n=1 \quad \text{and} \quad S^{a,v}_j=0 \quad \text{for} \quad j \ne n.$$
Now we claim that for every 0-1 solution $\left\{ x^{uv}_{ij}\right\}$  of the system (2.2.3)--(2.2.4), for any vertex $u \in V_1$, there is a unique vertex $i_u \in V_2$ 
such that the following equations hold:
\bigskip
\noindent (2.3.5) For all neighbors $v$ of $u$ we have 
$$S^{u,v}_{i_u}=1 \quad \text{and} \quad S^{u,v}_j=0 \quad \text{for} \quad j \ne i_u.$$
Then for the map $\phi: V_1 \longrightarrow V_2$ defined by $\phi(u)=i_u$ the conditions (2.3.1) are satisfied.

Clearly, if a choice $u \longmapsto i_u$ exists, it is unique. Because of (2.3.4), the equations (2.3.5) hold for $u=a$ and $i_u=n$. Since $G_1$ is connected, it suffices to show that whenever (2.3.5) holds for some vertex $u$ then for every neighbor $w$ of $u$ we can define $i_w \in V_2$ so that (2.3.5) holds with $u$ replaced by $w$ throughout.
Indeed, let $w$ be a neighbor of $u$ such that $(u, w) \in E_1$. It follows by (2.3.5) that there exists
$i_w$ such that 
$$x^{uw}_{i_u i_w} =1 \quad \text{and} \quad x^{uw}_{jk} =0 \quad \text{whenever} \quad j \ne i_u  \quad \text{or} \quad k \ne i_w.$$
From (2.3.3) it follows that for any neighbor $v$ of $w$, we have 
$$S^{w, v}_{i_w}=S^{w, u}_{i_w}=1 \quad \text{and} \quad S^{w, v}_{j} =S^{w, u}_j=0 \quad \text{for} \quad j \ne i_w, $$
as required. The case of neighbors $w$ of $u$ such that $(w, u) \in E_1$ is handled similarly. This proves that 0-1 solutions $\left\{x^{uv}_{ij}\right\}$, if 
any, of the system (2.3.3)--(2.3.4), are in one-to-one correspondence with graph homomorphisms $\phi: G_1 \longrightarrow G_2$ such that $\phi(a)=n$.

As we are interested in keeping the system (2.3.3)--(2.3.4) as sparse as possible, we arrange the equations (2.3.3) as follows: for a given $u \in V_1$, we list 
the neighbors $v$ of $u$ in some order $v_1, \ldots, v_m$ and then equate $S^{u, v_k}_i-S^{u, v_{k+1}}_i =0$ for $k=1, \ldots, m-1$. When the chosen vertex $a$ is a 
neighbor, we let $v_1=a$. This way the system (2.3.3)--(2.3.4) has not more than $2d_2$ variables per equation, where $d_2$ is the largest degree of a vertex of $G_2$, and not more than 4 equations per variable. 

Suppose that we are given a homomorphism $\phi: G_1 \longrightarrow G_2$ satisfying the constraint $\phi(a)=n$ for a fixed vertex $a$ of $G_1$ and a fixed vertex 
$n$ of $G_2$. 
As in Section 1.4, for an $\omega \in {\Bbb C}$ we consider the sum
$$\sum_{\psi:\ \psi(a)=n} \omega^{2\dist(\phi, \psi)}, \tag2.3.6$$
where $\psi$ ranges over all graph homomorphisms satisfying $\psi(a)=n$ and $\dist(\phi, \psi)$ is the number of directed edges where $\phi$ and $\psi$ disagree.
As follows from Section 1.4, we can approximate (2.3.6) within relative error $\epsilon >0$ in $d_2^{O(\ln |E_1| +\ln |E_2| -\ln \epsilon)}$
time provided
$$|\omega| \ \leq \ {\gamma \over d_2} \tag2.3.7$$
for some absolute constant $\gamma >0$ (we can choose $\gamma=0.1$). If the largest degree $d_2$ of a vertex of $G_2$ is fixed in advance, we obtain a
polynomial time approximation algorithm.

Suppose that $G_2$ is the complete graph with $n$ vertices and no loops, so that a homomorphism $G_1 \longrightarrow G_2$ is interpreted as a proper $n$-coloring of $G_1$ and $d_2=n-1$. If $n > d_1$, where $d_1$ is the largest degree of a vertex of $G_1$, it is trivial to come up with a homomorphism (proper $n$-coloring) $\phi: G_1 \longrightarrow G_2$ having a prescribed value on a prescribed vertex. In this case, the sum (2.3.6) is taken over all proper $n$-colorings 
$\psi$ of $G_2$ and each coloring is counted with weight exponentially small in the number of edges of $G_1$ whose coloring differ under $\phi$ and $\psi$.
If we could choose $\omega=1$ in (2.3.6), we would have counted all proper $n$-colorings of $G_1$ with $n > d_1$ colors, a notoriously difficult problem, see
\cite{Vi00} and \cite{C+19} for a randomized polynomial time approximation algorithm for counting $n$-colorings assuming that $n > (11/6) d_1$.

Given a pair of graphs $G_1$ and $G_2$, let us modify $G_2$ to a graph $\widehat{G}_2$ by adding an extra vertex $n+1$
with a loop and connected to all other vertices of $G_2$. Then there is always a homomorphism $\phi: G_1 \longrightarrow \widehat{G}_2$ 
which sends every vertex of $G_1$ to the newly added vertex $n+1$. In this case the sum (2.3.6) with $G_2$ replaced by $\widehat{G}_2$ and 
$n$ replaced by $n+1$ is interpreted as the sum over all homomorphisms of the induced subgraphs of $G_1$ to $G_2$. 

\subhead (2.4) Computing weight enumerators of linear codes \endsubhead If $\kappa$ is a prime, the set ${\Bbb Z}/\kappa {\Bbb Z}$ is identified with the finite field ${\Bbb F}_{\kappa}$ with $\kappa$ elements and $\left({\Bbb Z}/\kappa{\Bbb Z}\right)^n$ is the $n$-dimensional vector space over ${\Bbb F}_{\kappa}$. A set $X \subset {\Bbb F}_{\kappa}^n$ is called a 
 {\it code}.  The univariate polynomial
$$p_X(z)=1+\sum_{k=1}^n p_k(X) z^k,$$
where $p_k(X)$ is the number of vectors in $X$ with exactly $k$ non-zero coordinates, is called the {\it weight enumerator} of $X$, see for example, Chapter 3 of \cite{Li99}. 

Suppose that $X \subset {\Bbb F}_{\kappa}^n$ is defined by a system of linear equations
$$X=\left\{ x \in {\Bbb F}_{\kappa}^n: \quad Ax=0 \right\}, \tag2.4.1$$
where $A=\left(a_{ij}\right)$ is an $m \times n$ matrix with entries $a_{ij} \in {\Bbb F}_{\kappa}$. Hence $X \subset {\Bbb F}_{\kappa}^n$ is a subspace, called a {\it linear code}.
Generally, it is hard to compute $p_X(z)$ as it is hard to determine the smallest $k\geq 1$ with $p_k(X) \ne 0$, see \cite{B+78} and \cite{BN90}.

 Suppose now that the number of non-zero entries in every row of $A$ does not exceed $r\geq 2$ and the number of non-zero entries in every column of $A$ does not exceed $c \geq 1$. Let us define weights
$$w_1 = \ldots = w_n =z$$
for some $z \in {\Bbb C}$. Then 
$$w(X)=p_X(z)$$ and Theorem 1.6 implies that $p_X(z) \ne 0$ provided $|z| \leq \alpha/(\kappa-1)r \sqrt{c}$ and that 
$p_X(z)$ can be approximated within relative error $\epsilon >0$ in $(rc)^{O(\ln \kappa n - \ln \epsilon)}$ time, provided 
$|z| \leq \beta/(\kappa-1)r \sqrt{c}$, where $\beta < \alpha$ is fixed in advance. Again, if $r$ and $c$ are fixed in advance, we obtain an algorithm of polynomial $m(\kappa n/\epsilon)^{O(1)}$ complexity. Linear codes $X$ (typically binary, that is for $\kappa=2$) for which the number of non-zero entries in each row of the matrix $A$ in (2.4.1) is small are called {\it low-density parity-check} codes. They have many desirable properties and are of considerable interest, cf. Section 11 of \cite{MM09}. 

Let $C=X^{\bot}$, $C \subset {\Bbb F}_{\kappa}^n$, be the subspace (linear code) spanned by the rows of $A$ (we say that $A$ is the {\it generator matrix} of $C$). The MacWilliams identity  for the weight enumerators of $p_X$ and $p_C$ (see Theorem 3.5.3 of \cite{Li99}) states that 
$$p_X(z)={1 \over \kappa^{\dim C}} \bigl(1+(\kappa-1)z\bigr)^n p_C\left({1-z \over 1+(\kappa-1)z}\right).$$
It follows that 
$$p_C\left({1-z \over 1+(\kappa-1)z}\right) \ne 0 \quad \text{provided} \quad \left| z\right| \ \leq \ {\alpha \over (\kappa-1)r \sqrt{c}}$$
and that the value of 
$$p_C\left({1-z \over 1+(\kappa-1)z} \right)$$
 can be efficiently approximated provided 
 $$|z| \ \leq \ {\beta \over  (\kappa-1) r \sqrt{c}}.$$ In other words, the weight enumerator $p_C(z)$ of a linear code $C$ with a sparse code generator matrix is non-zero and can be efficiently approximated provided $|1-z|=O\left(1/r \sqrt{c}\right)$, where $r$ is an upper bound on the number of non-zero entries in every row, $c$ is an upper bound on the number of non-zero entries in every column of the matrix and the implied constant in the ``$O$" notation is absolute (in particular, it does not depend on $\kappa$). 
 
 One notable example of such a code with a sparse generating matrix is the binary {\it cut code} consisting of the indicators of cuts in a given graph $G=(V, E)$ with set $V$ of vertices and set $E$ of edges, see Section 1.9 of \cite{Di05} and \cite{BN90}, that is, indicators of subsets $E_S \subset E$ consisting of the edges with one endpoint in $S \subset V$ and the other in 
 $V \setminus S$. The rows of the code generating matrix are parameterized by vertices $v \in V$ of the graph, the columns are parameterized by the edges $e$ of the graph and the $(v, e)$ entry of the matrix is $1$ if $v$ is an endpoint of $e$ and 0 otherwise (hence each row is the indicator of the cut associated with the corresponding vertex).
We observe that the code generating matrix of a cut code contains at most $d(G)$ non-zero entries in every row, where $d(G)$ is the largest degree of a vertex of $G$, and exactly two non-zero entries in every column. The obtained algorithm for computing the weight of a cut code achieves roughly the same approximation as the algorithms of 
\cite{PR17a} and of Chapter 7 of \cite{Ba16}, where we approach computing weights of cuts via the graph homomorphism partition function.

\subhead (2.5) Ferromagnetic Potts model at low temperatures \endsubhead Let $G=(V, E)$ be a connected undirected graph, without loops or multiple edges. Given a real $\beta>0$ and an integer $\kappa > 1$, we consider the sum
$$P_{G, \kappa}(\beta)=\sum_{\phi: V \longrightarrow \{0, \ldots, \kappa-1\}} \exp\left\{ \beta \sum_{\{u, v\} \in E} \delta_{\phi(u) \phi(v)} \right\}, \tag2.5.1$$
where 
$$\delta_{ij}=\cases 1 &\text{if \ } i =j, \\ 0 &\text{if\ } i \ne j. \endcases$$
The expression (2.5.1) is known as the partition function of the {\it ferromagnetic} (since $\beta >0$) {\it Potts model with $\kappa$ colors}, see, for example, \cite{FV18}. Here the numbers $0, 1\ldots, \kappa-1$ 
are interpreted as colors: we color the vertices of $G$ with $\kappa$ colors in all possible ways, and each edge of $G$ with identically colored endpoints contributes to the inner sum. 
The number $\beta$ plays the role of the inverse temperature.
Using cluster expansions, it was shown in \cite{H+18} that for some induced subgraphs $G$ of the lattice ${\Bbb Z}^d$ the sum (2.5.1) can be approximated in polynomial time 
provided $\beta > \beta_0(d, \kappa)$ for some constant $\beta_0$ (that is, at sufficiently low temperatures). Here we deduce this result for a wide family of graphs and an explicit 
bound on $\beta_0$ from our Theorem 1.6. 

First, we rewrite (2.5.1) in the form
$$\split P_{G, \kappa}(\beta) = &e^{\beta |E|} \sum_{\phi: V \longrightarrow \{0, \ldots, \kappa-1\}} \prod_{\{u, v\} \in E} w(\phi(u), \phi(v)) \\
&\text{where} \quad w(i, j) = w_{\beta}(i, j)=e^{\beta(\delta_{ij}-1)}. \endsplit \tag2.5.2$$
Since $\beta >0$, we have $|w(i, j)| \leq 1$ and $w(i, j)=1$ if and only if $i=j$.

Next, we write the sum in (2.5.2) in the form $w(X)$, where $X$ is the set in Theorem 1.6. For that, we interpret colors $0, 1, \ldots, \kappa-1$ as remainders modulo $\kappa$.
We direct the edges of $G$ in an arbitrary way and with every, now directed, edge $(u, v)$ we associate a variable $x_{uv}$ taking values in ${\Bbb Z}/\kappa {\Bbb Z}$. 
The intended meaning of the variables $x_{uv}$ is that 
$$x_{uv} \equiv \phi(v) - \phi(u) \mod \kappa \quad \text{for all} \quad (u, v) \in E, \tag2.5.3$$
so that $x_{uv} \equiv 0$ if and only if the endpoints of the edge $\{u, v\}$ are colored with the same color. Given a set $\{x_{uv}: \ (u, v) \in E\}$, a solution 
$\phi: V \longrightarrow {\Bbb Z}/\kappa {\Bbb Z}$ to the system (2.5.3) exists, if and only if $\{x_{uv}\}$ satisfy the system of linear equations, constructed as follows:
we pick a cycle $C$ in $G$, orient it arbitrarily, and write
$$\sum\Sb \{u, v\} \in C: \\ (u, v) \text{\ is co-oriented with $C$} \endSb x_{uv} - \sum\Sb \{u, v\} \in C: \\ (u, v) \text{\ is counter-oriented with $C$} \endSb x_{uv} \equiv 0
\mod \kappa. \tag2.5.4$$ 
Moreover, since $G$ is connected, as long as the equations (2.5.4) are satisfied for all cycles $C$, the system (2.5.3) has exactly $\kappa$ solutions, that differ by a shift by an element of ${\Bbb Z}/\kappa {\Bbb Z}$. Indeed, if the equations (2.5.3) are satisfied then clearly (2.5.4) holds. On the other hand, given a solution to (2.5.4), we pick a vertex $v$ and assign the value of $\phi(v)$ arbitrarily. Then for every vertex $w$, we choose a path connecting $w$ to $v$, assign values of $\phi$ to the vertices along the path (in a necessarily unique way) so that the equations (2.5.3) are satisfied. Because of (2.5.4), the value of $\phi(w)$ does not depend on the chosen path.

Let $X \subset \left({\Bbb Z}/\kappa {\Bbb Z}\right)^{E}$ be the set of solutions of the system (2.5.4). We introduce a weight 
$w_{uv}=e^{-\beta}$
for each coordinate $x_{uv}$ with $(u, v) \in E$ and write (2.5.2) as 
$$P_{G,\kappa}(\beta)=\kappa e^{\beta |E|} w(X),$$
where $X$ is the set of solutions to the system (2.5.4). 

The equations (2.5.4) are not independent: it suffices to write (2.5.4) for a set of cycles ${\Cal C}$ that generate the homology group $H_1(G; {\Bbb Z})$. In view of Theorem 1.6, we would like to choose such a generating set ${\Cal C}$ of $H_1(G; {\Bbb Z})$ so that the number of edges in each cycle $C\in {\Cal C}$ does not exceed some $r \geq 2$ and the number of cycles $C\in {\Cal C}$ containing a given edge does not exceed some $c \geq 1$, for the smallest possible values of $r$ and $c$. Then we can approximate the partition function $P_{G, \kappa}(\beta)$ of (2.5.1) -- (2.5.2) provided 
$$\beta \ \geq \ 0.8+ \ln \left( (\kappa-1) r \sqrt{c}\right) \ >\ -\ln 0.45 + \ln \left( (\kappa-1) r \sqrt{c}\right),$$
and for fixed $r$ and $c$, we get a polynomial time approximation algorithm.

For example, suppose that $G$ is an induced subgraph of the integer lattice ${\Bbb Z}^d$ (with $d\geq 2$) constructed as follows. Given a point $(a_1, \ldots, a_d) \in {\Bbb Z}^d$, we call the set
$$\left\{ (x_1, \ldots, x_d): \ a_k\  \leq \ x_k \ \leq \ a_k+1: \ k=1, \ldots, d \right\}$$ 
 an {\it elementary cube}. We take finitely many elementary cubes whose union $U$ is a simply connected subset of ${\Bbb R}^d$ and let $G$ be the induced subgraph with vertices in $U$. Then there is a system of generators, ${\Cal C}$, of $H_1(G; {\Bbb Z})$ consisting of cycles with $r=4$ edges each and such that every edge of $C\in {\Cal C}$ belongs to at most $c=2(d-1)$ cycles (we choose the cycles on the boundary of 2-dimensional faces of the elementary cubes comprising $U$). Hence for such a graph $G$, we obtain a polynomial time approximation algorithm for $P_{G, \kappa}(\beta)$ provided
$$\beta \ \geq \ 2.6+  \ln \left( (\kappa-1) \sqrt{d-1}\right) \ > \ \ln \frac{4\sqrt{2}}{0.45} +  \ln \left( (\kappa-1) \sqrt{d-1}\right).$$

\head 3. Integrating over the torus \endhead

We begin our preparations to prove Theorems 1.2 and 1.6.
\subhead (3.1) Laurent polynomials on the torus \endsubhead
Let 
$${\Bbb S}^1 = \left\{z \in {\Bbb C}: \quad |z|=1 \right\}$$
be the unit circle in the complex plane and let 
$$\TT^m = {\Bbb S}^1 \times \ldots \times {\Bbb S}^1$$
be the direct product
of $m$ copies of ${\Bbb S}^1$ (torus), endowed with the product measure $\mu=\mu_1 \times \ldots \times \mu_m$, where $\mu_i$ is a Borel probability measure on the 
$i$-th copy of ${\Bbb S}^1$. We consider Laurent polynomials $p: \TT^m \longrightarrow {\Bbb C}$, 
$$p\left(z_1, \ldots, z_m \right)=\sum_{a \in A} \gamma_a \zz^{a} \tag3.1.1$$
as random variables on $\TT^m$. Here $A \subset {\Bbb Z}^m$ is a finite set of integer vectors, $\gamma_a \in {\Bbb C}$ for all $a \in A$ and
$$\zz^a=z_1^{\alpha_1} \cdots z_m^{\alpha_m} \quad \text{provided} \quad a=\left(\alpha_1, \ldots, \alpha_m \right),$$
where $z_i^0=1$. We are interested in conditions on the coefficients $\gamma_a$ which ensure that $\EE e^p \ne 0$. 

For $a \in A$ we define the {\it support} of $a$ by 
$$\supp a= \left\{i:\ \alpha_i \ne 0 \right\} \quad \text{where} \quad a=\left(\alpha_1, \ldots, \alpha_m\right).$$
Consequently, $|\supp a |$ is the number of non-zero coordinates of $a \in {\Bbb Z}^m$. In this section, we prove the following main result.
\proclaim{(3.2) Theorem} Let $p: \TT^m \longrightarrow {\Bbb C}$ be a Laurent polynomial as in (3.1.1).
Suppose that for some $0 \leq \theta_1, \ldots, \theta_m < 2\pi/3$, we have 
$$2 \sum\Sb a \in A: \\ i \in \supp a \endSb \left| \gamma_a \right| \prod_{j \in \supp a}{1 \over  \cos (\theta_j/2)}
 \ \leq \ \theta_i \quad \text{for} \quad i=1, \ldots, m.
\tag3.2.1$$
Then $$\EE e^p \ne 0.$$
\endproclaim

By choosing $\theta_i$ in a particular way, we obtain the following corollary.
\proclaim{(3.3) Corollary} There exists an absolute constant $\tau > 0$ such that if $p: \TT^m \longrightarrow {\Bbb C}$ is a Laurent polynomial as in (3.1.1)
and
$$|\supp a| \leq c \quad \text{for all} \quad  a \in A $$
and some $c \geq 1$ and 
$$\sum\Sb a \in A: \\  i \in \supp a \endSb \left| \gamma_a\right| \ \leq \ {\tau \over \sqrt{c}} \quad \text{for} \quad i=1, \ldots, m,$$
then 
$$\EE e^p \ne 0.$$
One can choose $\tau=0.56$.
\endproclaim

The proof is somewhat similar to that of \cite{Ba17} for $\EE e^p$ where $p: \{-1, 1\}^m \longrightarrow {\Bbb C}$ is a polynomial on the Boolean 
cube.

We start with a simple lemma (a discrete version of this lemma was suggested by Bukh \cite{Bu15}).
\proclaim{(3.4) Lemma} Let $f: \Omega \longrightarrow {\Bbb C}$ be a random variable and let $0 \leq \theta < 2\pi/3$ be a real number such that 
$f(\omega) \ne 0$ for all $\omega \in \Omega$ and the angle between any two complex numbers $f(\omega_1) \ne 0$ and $f(\omega_2) \ne 0$ considered as vectors in ${\Bbb R}^2 = {\Bbb C}$ does not exceed $\theta$. Suppose further that $\EE |f| < +\infty$. Then 
$$|\EE f| \ \geq \ \left(\cos {\theta \over 2} \right) \EE |f|.$$
\endproclaim 
\demo{Proof} First, we claim that $0$ does not lie in the convex hull of vectors $f(\omega) \in {\Bbb C}={\Bbb R}^2$. Otherwise we conclude by the Carath\'eodory Theorem that 
$0$ is a convex combination of some $3$ vectors $f(\omega_1)$, $f(\omega_2)$ and $f(\omega_3)$ and the angle between some two of them is at least $2\pi/3$, which is a contradiction. Hence the vectors $f(\omega)$ lie in some convex cone (angle) $K \subset {\Bbb C}$ measuring at most $\theta$ and with vertex at $0$. Let $\LL: {\Bbb R}^2 \longrightarrow {\Bbb R}^2$ be the orthogonal projection onto the bisector of $K$. Then 
$$\left| \EE f \right| \ \geq \ \left| \LL (\EE f )\right| = \left| \EE \LL(f) \right| =\EE |\LL(f)| \ \geq \EE \left( |f| \cos {\theta \over 2}\right) = 
\left( \cos {\theta \over 2}\right) \EE |f|.$$
Here the first (reading from left to right) inequality follows since the length of the orthogonal projection of a vector does not exceed the length of the vector; the next identity follows since $\LL$ is a linear operator; the next identity follows since for all $z \in K$ the vectors $\LL(z)$ are non-negative multiples of each other; the next inequality follows since 
$$\left| \LL(z) \right| \ \geq \left( \cos {\theta \over 2}\right) |z| \quad \text{for all} \quad z \in K;$$
 and the final identity follows since the expectation is a linear operator.
{\hfill \hfill \hfill} \qed
\enddemo

\subhead (3.5) Proof of Theorem 3.2 \endsubhead  For a function $f: \TT^m \longrightarrow {\Bbb C}$ and a subset $I \subset \{1, \ldots, m\}$, we denote by 
$\EE_I f$ the conditional expectation of $f$ obtained by integrating $f$ over the variables $z_i$ with $i \in I$. Hence if $f$ is a function of $z_1, \ldots, z_m$ and $I \subset \{1, \ldots, m\}$ then $h_I=\EE_I f$ is a function of $z_i$ for $i \notin I$. In particular, $h_I =f$ if $I=\emptyset$ and $h_I = \EE f$ if $I=\{1, \ldots, m\}$. 
If $I$ consists of a single element $i$, we write $\EE_i f$ instead of $\EE_{\{i\}} f$. We denote
$$\overline{I} = \{1, \ldots, m\} \setminus I$$
the complement of $I$. We will consider functions $f=e^p$ where $p: \TT^m \longrightarrow {\Bbb C}$ is a Laurent polynomial.

For $0 \leq \theta_1, \ldots, \theta_m < 2\pi/3$, we denote by $\PP_m\left(\theta_1, \ldots, \theta_m\right)$ the set of $m$-variate Laurent polynomials $p$ 
for which the inequalities (3.2.1) hold. Note that the condition $p \in \PP_m\left(\theta_1, \ldots, \theta_m\right)$ is a  finite system of linear inequalities for $|\gamma_a|$, $a \in A$.

Let us choose $p \in \PP_m\left(\theta_1,\ldots, \theta_m\right)$, let us fix some values $z_i \in {\Bbb S}^1$ for 
$I \subset \{1, \ldots, m\}$ and consider $p$ as a function of $z_i$ for $i\notin I$. It is not hard to see that  $p \in \PP_{m-|I|}\left(\theta_i: \ i \notin I\right)$.

We prove by induction on $m$ the following statements.
\bigskip
{\sl Statement $1_m$.} For any $p \in \PP_m\left(\theta_1, \ldots, \theta_m\right)$, we have $\EE e^p \ne 0$. Moreover, suppose that $p, q \in \PP_m \left(\theta_1, \ldots, \theta_m \right)$ are two Laurent polynomials that differ in at most one monomial, and the polynomial $p$ is obtained from $q$ by multiplying 
the coefficient $\gamma_b$ of some $\zz^b$ by some $\zeta \in {\Bbb S}^1$. Then the angle between $\EE e^p \ne 0$ and $\EE e^q \ne 0$ does not exceed
$$2 |\gamma_b| \prod_{i \in \supp b} {1 \over \cos (\theta_i/2)}.$$
\bigskip
{\sl Statement $2_m$.} Let $p \in \PP_m\left(\theta_1, \ldots, \theta_m \right)$ be a Laurent polynomial. Let $I=\{1, \ldots, m\} \setminus \{i\}$ for some 
$1 \leq i \leq m$ and let
$h_{I}(z_i) =\EE_I e^p$. Then for any $z_i', z_i'' \in {\Bbb S}^1$, we have $h_I(z_i') \ne 0$, $h_I(z_i'') \ne 0$ and the angle between the two complex numbers does not exceed $\theta_i$.
\bigskip
We start by proving Statement $2_1$. Then 
$$p(z) = \sum_{a \in A} \gamma_a z^a \quad \text{for some finite} \quad A \subset {\Bbb Z}$$ 
 is a univariate Laurent polynomial. For any $z \in {\Bbb S}^1$, we have 
$$\left| \arg e^{p(z)} \right| \ \leq \ \left| \Im\thinspace p(z) \right| \ \leq \ \left| p(z) \right| \ \leq \ \sum_{a \in A} |\gamma_a| \ \leq \ {1 \over 2}\theta_1$$
and the result is immediate.

Next, we prove that Statements $2_s$ for $s \leq m$ imply Statement $1_m$. 

\noindent Let us choose $p \in \PP_m\left(\theta_1, \ldots, \theta_m\right)$. For a set $I \subset \{1, \ldots, m\}$, let 
$$h_I\left(z_i:\ i \notin I \right) = \EE_I e^p.$$ 
Assuming that $I \ne \{1, \ldots, m\}$, let us pick an $i \notin I$. Then 
$$h_{I \cup \{i\}} = \EE_i h_I.$$
Let us fix variables $z_j \in {\Bbb S}^1$ with $j \notin I \cup \{i\}$ arbitrarily and consider $p$ as a Laurent polynomial from $\PP_r\left(\theta_k: \ k \in I \cup\{i\}\right)$ with 
$r=|I|+1$. Thus $h_I$ is a function of a single variable $z_i \in {\Bbb S^1}$ and by Statement $2_{r}$ for any two $z_i', z_i'' \in {\Bbb S}^1$, the angle between 
$h_I\left(z_i'\right) \ne 0$ and $h_I\left(z_i''\right) \ne 0$ does not exceed $\theta_i$. It follows from Lemma 3.4 that $h_{I \cup\{i\}}\left(z_j: \ j \notin I \cup\{i\} \right) \ne 0$
and, moreover, 
$$\left|h_{I\cup\{i\}}\right| =\left| \EE_i h_I\right| \ \geq \ \left( \cos {\theta_i \over 2}\right) \EE_i \left| h_I\right|>0.$$
Iterating, we obtain 
$$\left|h_{I \cup J}\right| =\left| \EE_J h_I \right| \ \geq \ \left(\prod_{j \in J} \cos  {\theta_j \over 2}\right) \EE _J|h_I|>0 \quad \text{provided} \quad J \cap I =\emptyset. \tag3.5.1$$
In particular, choosing $J=\overline{I}$, we obtain that $\EE e^p \ne 0$.

Suppose now that $p, q \in \PP_m\left(\theta_1, \ldots, \theta_m\right)$ where $p$ is obtained from $q$ by replacing a single monomial $\gamma_b \zz^b$ by 
$\gamma_b \zeta \zz^b$ for some $\zeta \in {\Bbb S}^1$. Let us fix all the remaining coefficients of $p$ and $q$ and consider $\EE e^p$ as a function of the coefficient $\gamma_b$ of $\zz^b$ as long as the resulting polynomial remains in $\PP_m\left(\theta_1, \ldots, \theta_m\right)$ (note that the set of admissible values of $|\gamma_b|$ is convex and includes $0$). Since $\EE e^p \ne 0$ for all $p \in \PP_m\left(\theta_1, \ldots, \theta_m\right)$, we can choose a continuous branch of 
$\ln \EE e^p$ as a function of $\gamma_b$. Then we have
$${\partial \over \partial \gamma_b} \ln \EE e^p = {{\partial \over \partial \gamma_b} \EE e^p \over \EE e^p}= {\EE\left(\zz^b e^p\right) \over \EE e^p}.$$
Let $I=\supp b$. Then 
$$\left| \EE \left(\zz^b e^p\right) \right| = \left| \EE_I \EE_{\overline{I}} \left(\zz^b e^p\right) \right| =\left| \EE_I \left(\zz^b \EE_{\overline{I}} e^p \right)\right| =
\left| \EE_I \zz^b h_{\overline{I}} \right| \ \leq \ \EE_I \left| h_{\overline{I}} \right|.$$
Similarly, 
$$\left|\EE e^p \right| = \left| \EE_I \EE_{\overline{I}} e^p \right| = \left| \EE_I h_{\overline{I}}\right| \ \geq \  \left(\prod_{i \in I} \cos {\theta_i \over 2} \right) 
\EE_I \left| h_{\overline{I}}\right|>0$$
by (3.5.1). Therefore,
$$\left| {\partial \over \partial \gamma_b} \ln \EE e^p \right| \ \leq \ \prod_{i \in \supp b} {1 \over \cos(\theta_i/2)}$$ 
and hence 
$$\left| \ln \EE e^p - \ln \EE e^q \right| \ \leq \ 2 |\gamma_b| \prod_{i \in \supp b} {1 \over \cos(\theta_i/2)}.$$
Statement $1_m$ now follows.

Next, we prove that Statement $1_m$ implies Statement $2_{m+1}$.

\noindent  Let $p \in \PP_{m+1}\left(\theta_1, \ldots, \theta_{m+1}\right)$ be a polynomial and let us choose 
an $1 \leq i \leq m+1$. Let $I =\{1, \ldots, m+1\} \setminus \{i\}$ and let $h_I(z_i) =\EE_I e^p$. If we fix $z_i$, we can consider $p$ as a Laurent polynomial 
in $\PP_m\left(\theta_j: \ j \ne i \right)$. Moreover, if we change the value of $z_i=z_i'$ to $z_i=z_i''$ then only the coefficients $\gamma_a$ of $p$ with $i \in \supp a$ are affected, and each of those coefficients gets multiplied by some $\zeta_a \in {\Bbb S}^1$. Repeatedly applying Statement $1_m$, we conclude that $h_I\left(z_i'\right) \ne 0$, 
$h_I\left(z_i''\right) \ne 0$ and the angle between the two complex numbers does not exceed
$$2 \sum\Sb a \in A: \\ i \in \supp a \endSb \left| \gamma_a\right| \prod_{j \in \supp a} {1 \over \cos(\theta_j/2)}, $$
which does not exceed $\theta_i$ by the definition of $\PP_{m+1}\left(\theta_1, \ldots, \theta_{m+1}\right)$, and Statement $2_{m+1}$ follows.

This concludes the induction and proves that $\EE e^p \ne 0$. 
{\hfill \hfill \hfill} \qed

\subhead (3.6) Proof of Corollary 3.3 \endsubhead Let us choose 
$$\theta_1=\ldots = \theta_m ={\delta \over \sqrt{c}}$$
for some $0 < \delta < 2 \pi/3$ to be determined later. To have the conditions of Theorem 3.2 satisfied, it suffices to have 
$$2\sum\Sb a \in A: \\ i \in \supp a \endSb  \left| \gamma_a \right| \left( \cos {\delta \over 2\sqrt{c}}\right)^{-c} \ \leq \ {\delta \over \sqrt{c}}
\quad \text{for} \quad i=1, \ldots, m.$$
Since
$$\left(\cos {\delta \over 2\sqrt{c}}\right)^{c} \ \geq \ \cos {\delta \over 2} \quad \text{for} \quad 0 \leq \delta \leq \pi,$$
see \cite{Ba17},
it suffices to have 
$$\sum\Sb a \in A: \\ i \in \supp a \endSb  |\gamma_a| \ \leq \ {\delta \cos (\delta/2) \over 2 \sqrt{c}} \quad \text{for} \quad i=1, \ldots, m.$$
Optimizing over $\delta$, we choose $\delta =1.72$ and 
$$\tau = {\delta \cos (\delta/2) \over 2} \approx 0.561,$$
which concludes the proof.
{\hfill \hfill \hfill} \qed
  
 \head 4. Proofs of Theorems 1.2 and 1.6  \endhead
 
 First, we prove Theorem 1.2. 
 
 Let $\TT^m={\Bbb S}^1 \times \ldots \times {\Bbb S}^1$ be the torus as in Section 3 and let us choose $\mu_i$ to be the rotation invariant (Haar) probability measure 
 on the $i$-th copy of ${\Bbb S}^1$. Let $\mu=\mu_1 \times \ldots \times \mu_m$ be the Haar probability measure on $\TT^m$.
 
 \proclaim{(4.1) Lemma} Let $X \subset \{0, 1\}^n$ be a set and $w_1, \ldots, w_n$ be complex weights as in Theorem 1.2. Let $a_j$, $j=1, \ldots, n$, be the columns of the matrix $A$, considered 
 as integer $m$-vectors and let us define a Laurent polynomial $q: \TT^m \longrightarrow {\Bbb C}$ by 
 $$q\left(z_1, \ldots, z_m \right) = \prod_{j=1}^n \left(1 + w_j \zz^{a_j}\right),$$
 where
 $$\zz^a=z_1^{\alpha_1} \cdots z_m^{\alpha_m} \quad \text{provided} \quad a=\left(\alpha_1, \ldots, \alpha_m\right).$$
 Then 
 $$w(X) = \EE q.$$
 \endproclaim
 \demo{Proof} Since for $a \in {\Bbb Z}^m$, we have
 $$\EE \zz^a =\cases 1& \text{if\ } a=0, \\ 0 &\text{if\ } a \ne 0, \endcases$$
expanding the product that defines $q$, we get
$$\EE q=\sum\Sb \xi_1, \ldots, \xi_n \in \{0, 1\}: \\  \xi_1 a_1 + \ldots + \xi_n a_n = 0 \endSb w_1^{\xi_1} \cdots w_n^{\xi_n} = w(X).$$
 {\hfill \hfill \hfill} \qed
 \enddemo
 
 \subhead (4.2) Proof of Theorem 1.2 \endsubhead
 
 \noindent  Let $q\left(z_1, \ldots, z_m \right)$ be the Laurent polynomial of Lemma 4.1, so that $w(X)= \EE q$. 
 Assuming that $|w_j| < 1$ for $j=1, \ldots, n$, we  write
$$ \ln q=\sum_{j=1}^n \ln \left(1 + w_j \zz^{a_j} \right) = \sum_{j=1}^n \sum_{k=1}^{\infty} (-1)^{k-1} {w_j^k  \zz^{k a_j} \over k}.$$
 For a positive integer $N$, let us define a Laurent polynomial
 $$p_N\left(z_1, \ldots, z_m\right) = \sum_{j=1}^n \sum_{k=1}^{N} (-1)^{k-1} { w_j^k \zz^{k a_j} \over k},$$
 which is just a truncation of the series expansion for $\ln q$. Let $\gamma_a \ne 0$ be the coefficient of the Laurent monomial $\zz^a$ in $p_N$.
Then
 $$\left| \supp a \right| \ \leq \ c$$
 and for $i=1, \ldots, m$, we have
 $$\sum_{a:\ i \in \supp a} \left| \gamma_a\right| \ \leq \ \sum_{j:\ a_{ij} \ne 0} \sum_{k=1}^N{ \left| w_j \right|^k \over k} \ \leq \ r  \max_{j=1, \ldots, n}
- \ln (1 -\left|w_j \right|)
 \ \leq \ {0.56 \over \sqrt{c}}$$
 as long as 
 $$|w_j| \ \leq \ {0.46 \over r \sqrt{c}} \quad \text{for} \quad j=1, \ldots, n. \tag4.2.1$$
 (We use that $-\ln (1-x) \leq 1.2x$ for $0 \leq x \leq 0.3$ and that $r \geq 2$.)
 Therefore, by Corollary 3.3, $\EE e^{p_N} \ne 0$ as long as (4.2.1) holds. On the other hand, $\EE e^{p_N}$ is an analytic function of $w_1, \ldots, w_n$ in the polydisc
 (4.2.1) and $\EE e^{p_N}$ converges to $\EE q$ uniformly on compact subsets of the polydisc. By the Hurwitz Theorem, see for example, Section 7.5 of \cite{Kr92}, we have either $\EE q \ne 0$ in the polydisc
 or $\EE q \equiv 0$ in the polydisc. Since for $w_1=\ldots =w_n=0$, we have $\EE q=1$, we conclude that $\EE q \ne 0$ provided (4.2.1) holds.
 {\hfill \hfill \hfill} \qed
 
 \subhead (4.3) Proof of Theorem 1.6 \endsubhead We modify the choice of the probability measure $\mu$ on $\TT^m$ as follows:
 we choose $\mu_i$ to be the uniform probability measure on the roots of unity of degree $\kappa$ and let $\mu=\mu_1 \times \ldots \times \mu_m$. 
 We note that for $a \in {\Bbb Z}^m$, $a=\left(\alpha_1, \ldots, \alpha_m\right)$, we have 
 $$\EE \zz^a = \cases 1 &\text{if\ } \alpha_i \equiv 0 \mod \kappa \quad \text{for} \quad i=1, \ldots, m, \\ 0 &\text{otherwise.} \endcases$$
 Given an $m \times n$ integer matrix $A=\left(a_{ij}\right)$, we define $q\left(z_1, \ldots, z_m\right)$ by 
 $$q\left(z_1, \ldots, z_m \right)=\prod_{j=1}^n \left(1 + w_j \zz^{a_j} + \ldots + w_j \zz^{(\kappa-1) a_j}\right).$$
Then
 $$\EE q=\sum\Sb \xi_1, \ldots, \xi_n: \\ \ a_{i1} \xi_1 + \ldots + \xi_n a_{in} \equiv 0 \mod \kappa \\ \text{for} \ i=1, \ldots, m, \\
   \xi_j \in \{0, 1, \ldots, \kappa-1\} \\ \text{for} \ j=1, \ldots, n \endSb  \prod_{j:\ \xi_j \ne 0} w_j = w(X).$$
  Assuming that $|w_j| < (\kappa-1)^{-1}$ for $j=1, \ldots, n$, we expand
  $$\split \ln q =&\sum_{j=1}^n \ln \left(1+ w_j \zz^{a_j} + \ldots + w_j \zz^{(\kappa-1) a_j}\right)\\=
 & \sum_{j=1}^n \sum_{s=1}^{\infty} (-1)^s {\left(w_j \zz^{a_j} + \ldots + w_j \zz^{(\kappa-1) a_j}\right)^s \over s}. \endsplit$$
 For a positive integer $N$, let us define a Laurent polynomial
 $$p_N\left(z_1, \ldots, z_m\right)=\sum_{j=1}^n \sum_{s=1}^N (-1)^s {\left(w_j \zz^{a_j} + \ldots + w_j \zz^{(\kappa-1) a_j}\right)^s \over s}. $$
 For every Laurent monomial $\zz^a$ which appears in $p_N$ with a coefficient $\gamma_a \ne 0$, we have $|\supp a| \leq c$. If $i \in \supp a$, then the coefficient of 
 $\zz^a$ in the polynomial $\left( \zz^{a_j} + \ldots +  \zz^{(\kappa-1) a_j}\right)^s$ is non-zero only if $a_{ij} \ne 0$. Hence
for $i=1, \ldots, m$, we
 have 
 $$\split \sum_{a: \ i \in \supp a} \left| \gamma_a \right| \ \leq \ &\sum_{j: \ a_{ij} \ne 0} \sum_{s=1}^N  {\bigl((\kappa -1) |w_j|\bigr)^s \over s} 
  \ \leq \ r \max_{j=1, \ldots, n} - \ln \bigl(1-(\kappa-1) |w_j|\bigr) \\ \leq \ & {0.56 \over \sqrt{c}} \endsplit$$
  provided
  $$\left| w_j\right| \ \leq \ {0.46 \over (\kappa-1) r \sqrt{c}} \quad \text{for} \quad j=1, \ldots, n.$$
 The proof is then concluded  as in Section 4.2.
 {\hfill \hfill \hfill} \qed
  
\head 5. Approximating $w(X)$ faster \endhead

Let $X \subset \{0, 1\}^n$ be the set defined in Theorem 1.2.  We assume that the $m \times n$ matrix $A$ has no zero rows or 
columns, see Section 1.3.  Recall that $r \geq 2$ is an upper bound on the number of non-zero entries in a row of $A$ and $c \geq 1$ is an upper bound on the number of non-zero entries in a column of $A$. As in Section 1.3, we define a univariate polynomial $w(X; \zeta)$, that is the weight of the set $X$ under the scaled 
weights $\zeta w_1, \ldots, \zeta w_n$, so $w(X; \zeta)$ is a polynomial of some degree $d \leq n$.
We let $f(\zeta)=\ln w(X; \zeta)$ for $\zeta$ in a neighborhood of $0$. 

Our goal is to show that the term $f^{(k)}(0)$ in the Taylor expansion (1.3.3) 
can be computed in $n (rc)^{O(k)}$ time, where we assume the standard RAM machine model with logarithmic-sized words, and additionally we assume that given a column index $j$ of the matrix $A=\left(a_{ij}\right)$ we can in time $O(c)$ compute the row indices $i$ such that $a_{ij}\neq 0$ (otherwise the running time is bounded by $nm(rc)^{O(k)}$).
We note that in this section, all the implied constants in the ``$O$" notation are absolute. 
In particular, if 
$k=O(\ln n-\ln \epsilon)$ as in Section 1.3, and $r$ and $c$ are fixed beforehand we obtain an algorithm of a polynomial in $n/\epsilon$ complexity.

Our algorithm heavily relies on the ideas of \cite{PR17a}, see also \cite{L+17}. 
\subhead (5.1) The idea of the algorithm \endsubhead  Since $w(X; 0)=1$, we can write
$$w(X; \zeta)=\prod_{i=1}^d \left(1 - {\zeta \over \zeta_i}\right),$$
where $\zeta_1, \ldots, \zeta_d \ne 0$ for some $d \leq n$ are the roots of $w(X; \zeta)$, listed with multiplicity. 
Then 
$$f(\zeta) = \sum_{i=1}^d  \ln \left(1-{\zeta \over \zeta_i}\right)$$ 
and 
$${f^{(k)}(0) \over k!}=-{1 \over k} \sum_{i=1}^d \zeta_i^{-k}.$$
We introduce the {\it power sums}
$$\sigma_k(A, w)= \zeta_1^{-k} + \ldots + \zeta_d^{-k}. \tag5.1.1$$
Hence our goal is to compute $\sigma_k(A, w)$ in $n (rc)^{O(k)}$ time.

The crucial feature of the power sums $\sigma_k(A, w)$ is that they are {\it additive functions} of $A$ as is explained below.

In what follows, we consider the set ${\Cal M}$ of integer matrices $A$ with rows and columns indexed by non-empty finite subsets of the set ${\Bbb N}$ of positive integers and without zero rows or columns. For non-empty finite subsets 
$R, C \subset {\Bbb N}$, an $R \times C$ integer-valued matrix $A \in {\Cal M}$ is a function $A: R \times C \longrightarrow {\Bbb Z}$ and we write the $(i, j)$-th 
entry of $A$ as $A(i, j)$ for $i \in R$ and $j \in C$. We fix complex weights $w_j$ and define 
$$ X_A=\Bigl\{ \bigl(\xi_j: \ j \in C\bigr) \in \{0, 1\}^C:  \quad \sum_{j \in C} A(i, j) \xi_j =0 \quad \text{for} \ i \in R \bigr\}. \tag5.1.2$$ 
Similarly, we define  univariate polynomials 
$$w\left(X_A; \zeta\right) = \sum\Sb x \in X_A \\ x=\left(\xi_j: \ j \in C\right) \endSb \prod_{j \in C} \left(\zeta w_j\right)^{\xi_j} \tag5.1.3$$
and define power sums $\sigma_k(A, w)$ by (5.1.1) where $\zeta_1, \ldots, \zeta_d$ are the roots of \newline $w\left(X_A; \zeta\right)$, listed with multiplicity.

Let $A_1, A_2 \in {\Cal M}$ be respectively $R_1 \times C_1$ and $R_2 \times C_2$ matrices. Suppose that $R_1 \cap R_2=\emptyset$ and $C_1 \cap C_2 =\emptyset$. We define the {\it direct sum} $A=A_1 \oplus A_2$ as the 
$R \times C$ matrix, where $R=R_1 \cup R_2$, $C=C_1 \cup C_2$ and 
$$A(i,j)=\cases A_1(i, j) &\text{if\ } i \in R_1 \ \text{and} \  j \in C_1, \\ A_2(i, j) &\text{if\ } i \in R_2 \ \text{and}\ j \in C_2, \\ 0 &\text{elsewhere.} \endcases$$
Clearly, $A \in {\Cal M}$. 

Let $A_1, A_2 \in {\Cal M}$ be matrices such that $A=A_1 \oplus A_2$ is defined. We observe that
$$w\left(X_A; \zeta\right)=w\left(X_{A_1}; \zeta\right) w\left(X_{A_2}; \zeta\right)$$
and hence
$$\sigma_k\left(A_1 \oplus A_2, w\right)= \sigma_k(A_1, w ) + \sigma_k(A_2, w). \tag5.1.4$$

 Given an $R \times C$ matrix $A \in {\Cal M}$ and an $R_1 \times C_1$ matrix $B \in {\Cal M}$, we define the {\it index}
$\ind(B, A) =1$ if $R_1 \subset R$, $C_1 \subset C$,
$$A(i, j) = B(i, j) \quad \text{for all} \quad i \in R_1 \quad \text{and all} \quad j \in C_1$$
and
$$A(i, j)=0 \quad \text{for all} \quad  i \in R \setminus R_1 \quad \text{and all} \quad j \in C_1.$$
Otherwise, we say that $\ind(B, A) =0$. 

We define a filtration 
$${\Cal M}_1 \subset {\Cal M}_2 \subset \ldots \subset {\Cal M}_k \subset  \ldots, $$
where ${\Cal M}_k \subset {\Cal M}$ consists of the matrices with at most $k$ columns.

In Lemma 5.3 below we show that we can write 
$$\sigma_k(A; w) = \sum_{B \in {\Cal M}_k} \ind(B, A) \mu_k(B, w) \quad \text{for all} \quad A \in {\Cal M} \tag5.1.5$$
and some complex numbers $\mu_k(B, w)$.
Although the sum in (5.1.5) contains infinitely many terms, for each $A \in {\Cal M}$, only finitely many terms are non-zero, so (5.1.5) is well-defined.

We say that a matrix $B \in {\Cal M}$ is {\it connected} if it cannot be represented as a direct sum $B = B_1 \oplus B_2$ 
for some matrices $B_1, B_2 \in {\Cal M}$ and {\it disconnected} otherwise. In Corollary 5.5  below, we deduce from the additivity property (5.1.4) that $\mu_k(B, w) =0$  in (5.1.5) unless $B$ is connected.
In Section 5.6 we show for any given $m \times n$ matrix $A$ with at most $r$ non-zero entries in each row and at most $c$ non-zero entries in each column the number of connected matrices $B \in {\Cal M}_k$ with $\ind(B, A) =1$ is at most $n (rc)^{O(k)}$ and that all such matrices $B$ can be found in $n(rc)^{O(k)}$ time. Finally, in Section 5.7 we show that for each connected $B \in {\Cal M}_k$, one can compute $\mu_k(B, w)$ in $cn 2^{O(k)}$ time. This produces an algorithm of $n(rc)^{O(k)}$ complexity for computing $\sigma_k(A, w)$. 
\bigskip
Next, we supply the necessary details. We start with a technical result describing how the function $\ind(B, \cdot)$ behaves under multiplication.
 Let $B_1 \in {\Cal M}$ be an $R_1 \times C_1$ matrix and let $B_2 \in {\Cal M}$ be an $R_2 \times C_2$ matrix. If the restrictions of $B_1$ and $B_2$ onto 
$(R_1 \cap R_2) \times (C_1 \cap C_2)$ coincide, we define the {\it connected sum} $B=B_1 \# B_2$, $B \in {\Cal M}$, as the $(R_1 \cup R_2) \times (C_1 \cup C_2)$ matrix such that 
$$B(i, j)=\cases B_1(i, j) &\text{if\ } i \in R_1 \text{\ and \ } j \in C_1, \\ B_2(i, j) &\text{if\ } i \in R_2  \text{\ and\ } j \in C_2, \\ 0 &\text{otherwise.} \endcases$$
In particular, if $R_1 \cap R_2 =\emptyset$ and $C_1 \cap C_2=\emptyset$ then $B_1 \# B_2 = B_1 \oplus B_2$ is the direct sum of $B_1$ and $B_2$.
\proclaim{(5.2) Lemma} Let $B_1 \in {\Cal M}$ be an $R_1 \times C_1$ matrix and let $B_2 \in {\Cal M}$ be an $R_2 \times C_2$ matrix.

Suppose that the following conditions (1) -- (3) are satisfied:
\roster
\item For all $i \in R_1 \cap R_2$ and all $j \in C_1 \cap C_2$ we have $B_1(i, j)=B_2(i, j)$;
\item For all $i \in R_1 \setminus R_2$ and all $j \in C_1 \cap C_2$ we have $B_1(i, j)=0$;
\item For all $i \in R_2 \setminus R_1$ and all $j \in C_1 \cap C_2$ we have 
$B_2(i, j)=0$.
\endroster
Then $B=B_1 \# B_2$ is defined and 
$$\ind(B_1, A) \ind(B_2, A)=\ind(B, A) \quad \text{for all} \quad A \in {\Cal M}.$$
If any of the conditions (1)--(3) is violated then 
$$\ind(B_1, A) \ind(B_2, A)=0 \quad \text{for all} \quad A \in {\Cal M}.$$
\endproclaim
\demo{Proof} Clearly, if (1) is violated then $\ind(B_1, A) \ind(B_2, A) =0$ for all $A \in {\Cal M}$. 

Suppose that (2) is violated. We assume that $R_1 \subset R$,  for $\ind(B_1, A)=0$ otherwise. If $\ind(B_2, A) =1$ then $A(i, j)=0$ for all $i \in R_1 \setminus R_2$ and all $j \in C_1 \cap C_2$ and hence $\ind(B_1, A) =0$
so that $\ind(B_1, A) \ind(B_2, A)=0$. Similarly, if (3) is violated then $\ind(B_1, A) \ind(B_2, A)=0$ for all $A \in {\Cal M}$.

Hence it remains to consider the case when (1)--(3) hold. Without loss of generality we assume that $R_1\cup R_2$ is a subset of the rows of $A$ and that $C_1\cup C_2$ is a subset of the columns of $A$.

If $\ind(B_1, A)=0$ for some $A \in {\Cal M}$ then either $B_1(i, j) \ne A(i, j)$ for some $i \in R_1$ and some $j \in C_1$ or $A(i,j) \ne 0$ for some $i \notin R_1$ and some 
$j \in C_1$. In either case $\ind(B, A)=0$. Similarly, if $\ind(B_2, A)=0$ then $\ind(B, A)=0$. If $\ind(B_1, A) =\ind(B_2, A)=1$ then $B(i,j)=A(i,j)$ for all 
$i \in R_1 \cup R_2$ and all $j \in C_1 \cup C_2$ while $A(i, j)=0$ for all $i \notin R_1 \cup R_2$ and all $j \in C_1 \cup C_2$ and hence $\ind(B, A)=1$
as well.
{\hfill \hfill \hfill} \qed
\enddemo

If the conditions (1)--(3) of Lemma 5.2 are satisfied, we say that the matrices $B_1$ and $B_2$ are {\it compatible} and denote it $B_1 \sim B_2$.
Now we are ready to prove the existence of a decomposition (5.1.5).

\proclaim{(5.3) Lemma} For a positive integer $k$ and a matrix $B \in {\Cal M}_k$ one can define complex numbers $\mu_k(B, w)$ 
so that (5.1.5) holds for all $A \in {\Cal M}$.
\endproclaim
\demo{Proof} 
We write the polynomial (5.1.3) in the monomial basis. Assuming that $A$ is an $R \times C$ matrix, we have
$$w\left(X_A; \zeta\right) =1 + \sum_{k=1}^n \pi_k(A, w) \zeta^k$$
where 
$$\pi_k(A, w)= \sum\Sb x=\left(\xi_j:\ j \in C\right): \\ x \in X_A, \\ \sum_{j \in C} \xi_j =k \endSb \prod_{j \in C} w_j^{\xi_j},$$
where $X_A$ is defined by (5.1.2). 

We say that a set $S \subset C$ is the {\it support} of a vector $x \in X_A$, 
$x=\left(\xi_j:\ j \in C \right)$, provided $\xi_j \ne 0$ if and only if $j \in S$. Clearly, the support of a vector $x$ contributing to $\pi_k(A, w)$ is a 
set $S \subset C$ satisfying $|S| \leq k$ and the vector $x_S=\left(\xi_j:\ j \in S\right)$ satisfies $A_S x_S =0$, where $A_S$ is the  $R \times S$ matrix consisting of the columns of $A$ with indices in $S$.

 This allows us to write
$$\pi_k(A, w) = \sum_{B \in {\Cal M}_k}  \ind(B, A) \lambda_k(B, w) \tag5.3.1$$
where for $R_1 \times C_1$ matrix $B$ we have
$$\lambda_k(B, w)= \sum\Sb x=\left(\xi_j: \ j \in C_1 \right):\\ x \in X_B, \\ \text{support of $x$ is $C_1$,} \\
\sum_{j \in C_1} \xi_j=k \endSb \prod_{j \in C_1} w_j^{\xi_j}. \tag5.3.2$$
 Although formally the sum (5.3.1) is infinite, for each $A \in {\Cal M}$ we have $\ind(B, A) \ne 0$ for only 
finitely many $B \in {\Cal M}$, so (5.3.1) is well-defined.

We observe that
$$\pi_k(A, w)=(-1)^k e_k\left( \zeta_1^{-1}, \ldots, \zeta_d^{-1}\right), $$
where $e_k$ is the $k$-th elementary symmetric function and $\zeta_1, \ldots, \zeta_d$ are the roots of $w\left(X_A; \zeta \right)$, listed with their multiplicities (recall that the constant term of $w\left(X_A; \zeta \right)$ is 1).
Therefore, the Newton identities imply that 
$$k \pi_k(A, w) = - \sum_{i=1}^k \pi_{k-i}(A, w)\sigma_i(A, w) \quad \text{for all} \quad k \geq 1, \tag5.3.3$$
where we define
$$\pi_0(A, w)=1.$$
We define 
$$\mu_1(B, w)=-\lambda_1(B, w) \quad \text{for} \quad B \in {\Cal M}_1.$$
 Assuming that $\mu_i(B, w)$ are defined for $B \in {\Cal M}_i$ and $i=1, \ldots, k-1$, for $k \geq 2$ we define for $B \in {\Cal M}_k$
$$\aligned \mu_{k}(B, w)= &-k \lambda_{k}(B, w) \\ &\quad - \sum\Sb B_1 \in {\Cal M}_{k-i}, B_2 \in {\Cal M}_i  \\ \text{\ for \ } 1 \leq i \leq k-1 : \\ B_1 \sim B_2 
\text{\ and } B_1 \# B_2 =B \endSb 
\lambda_{k-i}(B_1, w) \mu_i(B_2, w). \endaligned \tag5.3.4$$
Here the sum is taken over all distinct ordered pairs of compatible matrices $(B_1, B_2)$ such that  $B_1 \# B_2=B$ (in particular, we may have $B_1=B_2$).
We observe that for each $B$ the sum contains only finitely many terms, so $\mu_k(B, w)$ is well-defined. The identity (5.1.5) now follows from (5.3.1), (5.3.3) and Lemma 5.2.
{\hfill \hfill \hfill} \qed
\enddemo

Our next goal is to show that in (5.1.5) we have $\mu_k(B, w) \ne 0$ only for connected matrices $B$. We start with a general structural result, very similar in spirit to 
Lemma 4.2 of \cite{CS16}, see also \cite{PR17a}.
\proclaim{(5.4) Lemma} Let us consider a function $f: {\Cal M} \longrightarrow {\Bbb C}$ defined by
$$f(A)=\sum_{B \in {\Cal S}} \mu_B \ind(B, A),$$
where ${\Cal S} \subset {\Cal M}$ is a (possibly infinite) set and $\mu_B \in {\Bbb C} \setminus \{0\}$ for all $B \in {\Cal S}$ (for each $A \in {\Cal M}$ only finitely many summands are non-zero, so $f$ is well-defined).
Suppose that 
$$f(A_1 \oplus A_2) =f(A_1) + f(A_2)$$
for any two matrices $A_1, A_2 \in {\Cal M}$ such that $A_1 \oplus A_2$ is defined. Then each $B \in {\Cal S}$ is connected, that is, cannot be written 
as $B=B_1 \oplus B_2$ for some $B_1, B_2 \in {\Cal M}$.
\endproclaim
\demo{Proof} 
Seeking a contradiction, assume that there is a disconnected $B \in {\Cal S}$.
We observe that if $B \in {\Cal M}$ is connected, then 
$$\ind(B, A_1 \oplus A_2) =\ind(B, A_1) + \ind(B, A_2)$$
for any two $A_1, A_2 \in {\Cal M}$ such that $A_1 \oplus A_2$ is defined (since $B$ is connected, we cannot have $\ind(B, A_1)=\ind(B, A_2)=1$ 
provided $A_1 \oplus A_2$ is defined). Therefore, without loss of generality, we assume that all $B \in {\Cal S}$ are disconnected. 
Let 
us choose a $D \in {\Cal S}$ that has the smallest number of columns. Hence we have 
$D=D_1 \oplus D_2$ for some $D_1$ and $D_2$. Then 
$$f(D_i)=\sum_{B \in {\Cal S}} \mu_B \ind(B, D_i)=0 \quad \text{for} \quad i=1, 2$$
since $D_i$ has fewer columns than any matrix $B \in {\Cal S}$. Therefore,
$$f(D) = f(D_1) + f(D_2) =0.$$
On the other hand, $\ind(B, D) =0$ for all $B \in {\Cal S}\setminus \{D\}$ and 
$$f(D) = \mu_D \ind(D, D) =\mu_D \ne 0,$$
which is a contradiction.
{\hfill \hfill \hfill} \qed
\enddemo

\proclaim{(5.5) Corollary} In the expansion (5.1.5), we have 
$\mu_k(B, w)=0$ whenever $B$ is disconnected. 
\endproclaim
\demo{Proof} Follows by (5.1.4) and Lemma 5.4.
{\hfill \hfill \hfill} \qed
\enddemo

\subhead (5.6) Enumerating connected matrices \endsubhead  Given an integer $k \geq 1$ and an $R \times C$ matrix $A \in {\Cal M}$ with at most $r$ non-zero entries in each row and at most $c$ non-zero entries in each column, we want to compile a list of all connected 
matrices $B \in {\Cal M}_k$ such that $\ind(B, A)=1$. First, we observe that an $R_1 \times C_1$ matrix 
$B \in {\Cal M}$ such that $\ind(B, A)=1$ is uniquely determined by its set of columns $C_1 \subset C$ since $R_1 \subset R$ is then the set of rows of $A$ whose 
restriction onto $C_1$ are not zero. 

We define a graph $G=(C, E)$. The vertices of $G$ are the columns of $A$ and two vertices $c_1$ and $c_2$ span an edge of $G$ if and only if there is a row of $A$ with non-zero entries in columns $c_1$ and $c_2$. We note that the degree of each vertex of $G$ does not exceed $d=r c$. To enumerate connected matrices $B \in {\Cal M}_k$ such that $\ind(B, A)=1$ is to enumerate sets of vertices of cardinality at most $k$ in $C$ that induce a connected subgraph of $G$. This is done as in \cite{PR17b}. The crucial observation is that as long as one vertex $c$ is 
chosen, there are at most 
$$k^{-1} {k d \choose k-1} \ \leq \ {(e d)^{k-1} \over 2},$$
 connected induced subgraphs with $k\geq 2$ vertices containing $c$, see Lemma 2.1 of \cite{B+13}. Consequently, there are  $n d^{O(k)}$ induced connected subgraphs with at most $k$ vertices in $G$.
 Once the vertex $c$ is chosen, the subgraphs are enumerated with $d^{O(k)}$ complexity, by successively exploring adjacent vertices, see \cite{PR17b} for details.

\subhead (5.7) Summary of the algorithm \endsubhead Given an $m \times n$ matrix $A$ without zero rows and columns, we interpret it as an $R \times C$ matrix $A \in {\Cal M}$, where $R=\{1, \ldots, m\}$ and $C=\{1, \ldots, n\}$. Given a positive integer $k$, as in Section 5.6 we compile a list ${\Cal C}$ of all connected matrices $B \in {\Cal M}_k$ such that $\ind(B, A)=1$. We define the filtration 
$${\Cal C}_1 \subset {\Cal C}_2 \subset \ldots \subset {\Cal C}_{k-1} \subset {\Cal C}_k={\Cal C},$$
where ${\Cal C}_i$ is the set of matrices $B \in {\Cal C}$ with at most $i$ columns.

 Given complex numbers $w_1, \ldots, w_n$, from Lemma 5.3 and (5.3.2) in particular, we obtain
 $$\mu_1(B, w)=-\lambda_1(B, w)=0 \quad \text{for} \quad B \in {\Cal C}_1,$$
since by our assumption $B$ has no zero rows.

Suppose that we have computed $\mu_i(B, w)$ for $i=1, \ldots, k-1$ and all $B \in {\Cal C}_{k-1}$ for $k \geq 2$.
To compute $\mu_k(B, w)$ for all $B \in {\Cal C}_k$,  we use formula (5.3.4). Since every matrix $B \in {\Cal C}_k$ has at most $k$ columns, 
there are not more than $4^k$ pairs of matrices $B_1\in {\Cal M}_{k-i}$ and $B_2 \in {\Cal M}_i$  such that $B_1 \# B_2 =B$ and all such pairs can be found by inspection in $O(4^k)$ time.
We then use (5.3.2) to compute the terms $\lambda_{k-i}(B_1, w)$ for $i=0, \ldots, k-1$. We note that there are ${k-i-1 \choose |S|-1} \leq 2^{k-i}$ non-negative integer vectors with support $S$ and the sum $k-i$ of the coordinates, so each $\lambda_{k-i}(B_1, w)$ is computed in $c(k-i)n2^{O(k-i)}$ time.

This gives us the list of values $\mu_k(B, w)$ for all $B \in {\Cal C}_k$. We then compute 
$$\sigma_k(A, w) =\sum_{B \in {\Cal C}_k} \mu_k(B, w),$$
as desired.

\head Acknowledgments \endhead 

The authors are grateful to Martin Dyer for pointing out to possible connections with \cite{Va79} and \cite{VV86}, to Alex Samorodnitsky for pointing out to possible connections with low-density parity-check codes, to Prasad Tetali for pointing us to \cite{BS94} and to Matthew Jenssen for telling us how to improve the degree dependence from $d_2-d_1\geq 2$ to $d_2-d_1\geq 1$ in Section 2.6. We thank the authors of \cite{C+16} for pointing out that their construction implies that computing the partition function in the hard-core model at high fugacity is a $\#$BIS-hard problem in the class of regular bipartite graphs. We thank the anonymous referee for catching an inaccuracy in one of our estimates.

\enddocument

\end